\documentclass[reqno]{amsart}
\usepackage{amssymb,url}
\usepackage[colorlinks=true,hypertex]{hyperref}
\date{This manuscript was completed on 7 July 2007}
\date{}
\allowdisplaybreaks[4]
\theoremstyle{plain}
\newtheorem{thm}{Theorem}
\newtheorem{lem}{Lemma}
\newtheorem{cor}{Corollary}
\theoremstyle{remark}
\newtheorem{rem}{Remark}
\DeclareMathOperator{\td}{d\mspace{-1mu}}
\newcommand{\tn}{\mathbb{N}}

\numberwithin{equation}{section}

\begin{document}

\title{Some properties of the psi and polygamma functions}

\author[F. Qi]{Feng Qi}
\address[F. Qi]{Research Institute of Mathematical Inequality Theory\\ Henan Polytechnic University\\ Jiaozuo City, Henan Province, 454010\\ China}
\email{\href{mailto: F. Qi <qifeng618@gmail.com>}{qifeng618@gmail.com}, \href{mailto: F. Qi <qifeng618@hotmail.com>}{qifeng618@hotmail.com}, \href{mailto: F. Qi <qifeng618@qq.com>}{qifeng618@qq.com}}
\urladdr{\url{http://qifeng618.spaces.live.com}}

\author[B.-N. Guo]{Bai-Ni Guo}
\address[B.-N. Guo]{School of Mathematics and Informatics\\
Henan Polytechnic University\\ Jiaozuo City, Henan Province, 454010\\ China} \email{\href{mailto: B.-N.
Guo <bai.ni.guo@gmail.com>}{bai.ni.guo@gmail.com}, \href{mailto: B.-N. Guo
<bai.ni.guo@hotmail.com>}{bai.ni.guo@hotmail.com}} \urladdr{\url{http://guobaini.spaces.live.com}}

\begin{abstract}
In this paper, some monotonicity and concavity results of several functions involving the psi and polygamma functions are proved, and then some known inequalities are extended and generalized.
\end{abstract}



\thanks{The first author was supported in part by the China Scholarship Council and the authors were supported in part by the NSF of Henan University, China}

\thanks{This paper was typeset using \AmS-\LaTeX}

\maketitle

\section{Introduction}

It is well-known that the classical Euler's gamma function $\Gamma(x)$ plays a central role in the theory of special functions and has much extensive applications in many branches, for example, statistics, physics, engineering, and other mathematical sciences. The logarithmic derivative of $\Gamma(x)$, denoted by
$\psi(x)=\frac{\Gamma'(x)}{\Gamma(x)}$, is called the psi or digamma function, and $\psi^{(i)}(x)$ for $i\in\mathbb{N}$ are known as the polygamma or multigamma functions.
\par
In \cite[Theorem~2]{batir-new}, it was discovered that if $a\le-\ln2$ and $b\ge0$, then
\begin{equation}\label{batir-ineq-orig}
a-\ln(e^{1/x}-1)<\psi(x)<b-\ln(e^{1/x}-1)
\end{equation}
holds for $x>0$. In \cite[Theorem~2.8]{batir-jmaa-06-05-065}, inequality \eqref{batir-ineq-orig} was sharpened: If $a\le-\gamma$ and $b\ge0$, then inequality \eqref{batir-ineq-orig} is valid for $x>0$, where the constants $-\gamma=-0.577\dotsc$, the negative of Euler-Mascheroni's constant, and $0$ are the best possible.
\par
The first aim of this paper is to generalize inequality \eqref{batir-ineq-orig} to an increasingly monotonic and concave properties as follows.

\begin{thm}\label{batir-thm2.8-gen}
The function
\begin{equation}
\phi(x)=\psi(x)+\ln\bigl(e^{1/x}-1\bigr)
\end{equation}
is not only strictly increasing but also strictly concave on $(0,\infty)$, with
\begin{equation}
\lim_{x\to0^+}\phi(x)=-\gamma\quad\text{and}\quad \lim_{x\to\infty}\phi(x)=0.
\end{equation}
\end{thm}

\begin{rem}
It is noted that the increasingly monotonic property in Theorem~\ref{batir-thm2.8-gen} was also obtained in \cite{alzer-expo-math-2006} by using a different approach.
\end{rem}

As direct consequences of the proof of Theorem~\ref{batir-thm2.8-gen}, the following two inequalities for the trigamma function $\psi'(x)$ and the tetragamma function $\psi''(x)$ are deduced.

\begin{cor}\label{batir-thm2.8-gen-cor}
For $x>0$,
\begin{gather}
\psi'(x)>\frac{e^{1/x}}{\bigl(e^{1/x}-1\bigr)x^2} \quad\text{and}\quad
\psi''(x)<\frac{e^{1/x}\bigl[1-2x\bigl(e^{1/x}-1\bigr)\bigr]} {\bigl(e^{1/x}-1\bigr)^2x^4}.
\end{gather}
\end{cor}

The second aim of this paper is to extend Theorem~\ref{batir-thm2.8-gen} to the following necessary and sufficient conditions.

\begin{thm}\label{batir-thm2.8-gen-2}
For $\theta>0$, let
\begin{equation}
\phi_\theta(x)=\psi(x)+\ln\bigl(e^{\theta/x}-1\bigr)
\end{equation}
on $(0,\infty)$.
\begin{enumerate}
\item
The function $\phi_\theta(x)$ is strictly increasing if and only if $0<\theta\le1$ and strictly decreasing if $\theta\ge2$.
\item
The function $\phi_\theta(x)$ is strictly concave if $0<\theta\le1$ and strictly convex if $\theta\ge2$.
\item
$\lim_{x\to\infty}\phi_\theta(x)=\ln\theta$ and
\begin{equation}
\lim_{x\to0^+}\phi_\theta(x)=\begin{cases}
-\gamma,&\theta=1,\\
\infty,&\theta>1,\\
-\infty,&0<\theta<1.
\end{cases}
\end{equation}
\end{enumerate}
\end{thm}

As straightforward consequences of the proof of Theorem~\ref{batir-thm2.8-gen-2}, the following inequalities for the trigamma function $\psi'(x)$ and the tetragamma function $\psi''(x)$ are presented, which extend the two inequalities in Corollary~\ref{batir-thm2.8-gen-cor}.

\begin{cor}\label{batir-thm2.8-gen-cor-2}
For $x>0$, inequalities
\begin{gather}\label{phi'-c}
\psi'(x)>\frac{\theta e^{\theta/x}} {x^2\bigl(e^{\theta/x}-1\bigr)} \quad\text{and}\quad
\psi''(x)\le\frac{\theta e^{\theta/x} \bigl[\theta-2x\bigl(e^{\theta/x}-1\bigr)\bigr]} {x^4\bigl(e^{\theta/x}-1\bigr)^2}
\end{gather}
hold if $0<\theta\le1$ and reverse if $\theta\ge2$.
\end{cor}
\par
In \cite[Theorem~2.6]{batir-jmaa-06-05-065}, inequality
\begin{equation}\label{batir-error-ineq}
-\gamma+x\psi'\biggl(\frac{x}2\biggr)<\psi(x+1) <-\gamma+x\psi'\bigl(\sqrt{x+1}\,-1\bigr)
\end{equation}
for $x>0$ was showed. Careful observation reveals that inequality \eqref{batir-error-ineq} is not valid: If taking $x=1$, the left-hand side inequality in \eqref{batir-error-ineq} is reduced to $-\gamma+\frac{\pi^2}2<1-\gamma$ which does not hold true clearly. After checking up the proof of \cite[Theorem~2.6]{batir-jmaa-06-05-065}, it is found that inequality \eqref{batir-error-ineq} can be corrected and extended as the following theorem.

\begin{thm}\label{batir-thm-corr}
If $x>0$,
\begin{equation}\label{batir-corr-ineq}
-\gamma+x\psi'\biggl(1+\frac{x}2\biggr)<\psi(x+1) <-\gamma+x\psi'\bigl(\sqrt{x+1}\,\bigr).
\end{equation}
If $-1<x<0$, inequality \eqref{batir-corr-ineq} is reversed.
\end{thm}

As a generalization of Theorem~\ref{batir-thm-corr}, the following monotonicity are obtained.

\begin{thm}\label{f-g-mon}
The functions
\begin{equation}\label{f-g-dfn}
f(x)=\psi(x+1)-x\psi'\biggl(1+\frac{x}2\biggr)\quad\text{and}\quad g(x)=x\psi'\bigl(\sqrt{x+1}\,\bigr)-\psi(x+1)
\end{equation}
are both strictly increasing on $(-1,\infty)$, with limits
\begin{align*}
\lim_{x\to-1^+}f(x)&=-\infty, & \lim_{x\to-1^+}g(x)&=1+\gamma-\frac{\pi^2}6,& \lim_{x\to\infty}f(x)&=\lim_{x\to\infty}g(x)=\infty.
\end{align*}
\end{thm}

\begin{rem}
Making use of the difference equation \eqref{psisymp4} below, the functions $f(x)$ and $g(x)$ defined in \eqref{f-g-dfn} for $x\in(-1,\infty)$ can be rewritten as
\begin{equation}
f(x)=\begin{cases}
\psi(x)-x\psi'\biggl(\dfrac{x}2\biggr)+\dfrac5{x},&x\ne0\\-\gamma,&x=0
\end{cases}
\end{equation}
and
\begin{equation}
g(x)=\begin{cases}
x\psi'\bigl(\sqrt{x+1}\,\bigr)-\psi(x)-\dfrac1x,&x\ne0\\
\gamma,&x=0.
\end{cases}
\end{equation}
\end{rem}

As an immediate consequence of the proof of Theorem~\ref{f-g-mon}, the monotonicity of the function
$$
\bigl(u^2-1\bigr)\psi'(u)-\psi\bigl(u^2\bigr)
$$
on $(-1,\infty)$ is derived as follows.

\begin{thm}\label{correction-thm-3}
For $x>-1$, the function
\begin{equation}
h(x)=\begin{cases}\bigl(x^2-1\bigr)\psi'(x)-\psi\bigl(x^2\bigr),&x\ne0\\
1+\gamma-\dfrac{\pi^2}6,&x=0
\end{cases}
\end{equation}
is strictly increasing on $(-1,\infty)$, with
\begin{equation}
\lim_{x\to-1^+}h(x)=-\infty\quad \text{and}\quad \lim_{x\to\infty}h(x)=\infty.
\end{equation}
\end{thm}

\begin{rem}
It is conjectured that the function $h(x)$ is strictly concave on $(-1,1)$ and strictly convex on $(1,\infty)$.
\end{rem}

Finally, Theorem~\ref{f-g-mon} can be generalized as follows.

\begin{thm}\label{f-g-mon-poly}
If $i$ is a positive odd, then the function
\begin{equation}\label{f-i-dfn}
f_i(x)=\psi^{(i)}(x+1)-x\psi^{(i+1)}\biggl(1+\frac{x}2\biggr)
\end{equation}
is strictly decreasing on $(-1,\infty)$; If $i$ is a positive even, then the function $f_i(x)$ is strictly increasing on $(-1,\infty)$; For all $i\in\mathbb{N}$, the limits
\begin{equation}
\begin{aligned}
\lim_{x\to-1^+}f_i(x)&=(-1)^{i+1}\infty&\text{and}&& \lim_{x\to\infty}f_i(x)&=0
\end{aligned}
\end{equation}
hold true.
\end{thm}

\begin{rem}
Similar to the monotonic properties of the function $f_i(x)$, the following conjecture can be posed: If $i$ is a positive odd, then the function
\begin{equation}\label{g-i-dfn}
g_i(x)=\psi^{(i)}(x+1)-x\psi^{(i+1)}\bigl(\sqrt{x+1}\,\bigr)
\end{equation}
is strictly decreasing on $(-1,\infty)$; If $i$ is a positive even, then the function $g_i(x)$ is strictly increasing on $(-1,\infty)$; For all $i\in\mathbb{N}$, the limits
\begin{equation}
\begin{aligned}
\lim_{x\to-1^+}g_i(x)&=(-1)^{i+1}\infty&\text{and}&& \lim_{x\to\infty}g_i(x)&=0
\end{aligned}
\end{equation}
are valid, except $\lim_{x\to\infty}g_1(x)=1$.
\par
Direct calculation yields
\begin{align*}
[g_i(x)]'&=\psi^{(i+1)}(x+1) -\psi^{(i+1)}\bigl(\sqrt{x+1}\,\bigr) -\frac{x}{2\sqrt{x+1}\,}\psi^{(i+2)}\bigl(\sqrt{x+1}\,\bigr)\\
&=\frac{2u\bigl[\psi^{(i+1)}\bigl(u^2\bigr) -\psi^{(i+1)}(u)\bigr] -\bigl(u^2-1\bigr)\psi^{(i+2)}(u)}{2u}\\
&=\frac{\bigl[\psi^{(i)}\bigl(u^2\bigr) -\bigl(u^2-1\bigr)\psi^{(i+1)}(u)\bigr]'}{2u},
\end{align*}
where $u=\sqrt{x+1}\,>0$ for $x>-1$. Therefore, in order to verify above conjecture, it is sufficient to show the monotonic properties of the function
\begin{equation}
\psi^{(i)}\bigl(u^2\bigr) -\bigl(u^2-1\bigr)\psi^{(i+1)}(u)
\end{equation}
on $(0,\infty)$.
\end{rem}

\begin{rem}
It is also natural to pose an open problem: For $i,k\in\mathbb{N}$ and positive numbers $\alpha$, $\beta$, $\delta$, $\lambda$, $\mu$ and $\tau$, what about the monotonicities and convexities of the more general function
\begin{equation}
\varphi_{i,k}(x)=\psi^{(i-1)}(x+\alpha) -(x+\beta)^k\psi^{(i)}\bigl(\lambda(x+\delta)^\mu+\tau\bigr)
\end{equation}
in an appropriate interval where it is defined?
\end{rem}

\section{Lemmas}

The following lemmas are useful for the proofs of some of our theorems.

\begin{lem}[{\cite[Lemma~1]{notes-best-simple-equiv.tex}}]\label{key-lem}
If $f(x)$ is a function defined in an infinite interval $I$ such that
$$
f(x)-f(x+\varepsilon)>0\quad \text{and}\quad \lim_{x\to\infty}f(x)=\delta
$$
for some $\varepsilon>0$, then $f(x)>\delta$ on $I$.
\end{lem}

\begin{proof}
By induction, for any $x\in I$,
\begin{equation*}
f(x)>f(x+\varepsilon)>f(x+2\varepsilon)>\dotsm>f(x+k\varepsilon)\to\delta
\end{equation*}
as $k\to\infty$. The proof of Lemma \ref{key-lem} is complete.
\end{proof}

\begin{lem}[\cite{abram}]
For $x>0$ and $k\in\tn$,
\begin{gather}
\label{binetcor}
\psi(x)-\ln x+\frac1x =\int_{0}^{\infty} \left(\frac{1}{t}-\frac{1}{e^t-1}\right)e^{-xt}\td t,\\
\label{psisymp4}
\psi^{(k-1)}(x+1)=\psi^{(k-1)}(x)+\frac{(-1)^{k-1}(k-1)!}{x^k}.
\end{gather}
\end{lem}

Recall \cite{clark-ismail,widder} that a function $f$ is said to be completely monotonic on an interval $I$ if it has derivatives of all orders on $I$ and
\begin{equation}
(-1)^{n}f^{(n)}(x)\ge0
\end{equation}
for $x \in I$ and $n \geq0$. The well-known Bernstein's Theorem \cite[p.~161]{widder} states that a function $f$ is completely monotonic on $(0,\infty)$ if and only if
\begin{equation}
f(x)=\int_0^\infty e^{-xs}\td\mu(s),
\end{equation}
where $\mu$ is a nonnegative measure on $[0,\infty)$ such that the integral converges for all $x>0$. This expresses that a function $f$ is completely monotonic on $(0,\infty)$ if and only if it is a Laplace transform of the measure $\mu$.

\begin{lem}[{\cite[Theorem~2]{sandor-gamma-2}}]\label{comp-thm-1}
The function
$$
\psi(x)-\ln x+\frac{\alpha}x
$$
is completely monotonic on $(0,\infty)$ if and only if $\alpha\ge1$ and
$$
\ln x-\frac{\alpha}x-\psi(x)
$$
is completely monotonic on $(0,\infty)$ if and only if $\alpha\le\frac12$. Consequently, inequalities
\begin{equation}\label{psi'ineq}
\frac{(k-1)!}{x^k}+\frac{k!}{2x^{k+1}}<
(-1)^{k+1}\psi^{(k)}(x)<\frac{(k-1)!}{x^k}+\frac{k!}{x^{k+1}}
\end{equation}
hold for $x\in(0,\infty)$ and $k\in\tn$.
\end{lem}

\begin{rem}
Recall \cite{Atanassov, CBerg, grin-ismail, compmon2, minus-one, e-gam-rat-comp-mon, auscmrgmia} that a positive function $f$ on an interval $I$ is called logarithmically completely monotonic if it satisfies
\begin{equation}
(-1)^k[\ln f(x)]^{(k)}\ge0
\end{equation}
for $k\in\mathbb{N}$ on $I$. Lemma~\ref{comp-thm-1} can also be concluded from the necessary and sufficient conditions such that the function
\begin{equation}
\frac{e^x\Gamma(x)}{x^{x-\alpha}}
\end{equation}
for $\alpha\in\mathbb{R}$ and its reciprocal are logarithmically completely monotonic on $(0,\infty)$, which was established in \cite{chen-qi-log-jmaa}.
\end{rem}

\begin{lem}[{\cite[pp.~259\nobreakdash--260]{abram}}]
For $z\ne-1,-2,-3,\dotsc$,
\begin{equation}\label{psi-sum}
\psi(1+z)=-\gamma+\sum_{n=1}^{\infty}\frac{z}{n(n+z)}.
\end{equation}
For $n\in\tn$ and $z\ne0,-1,-2,\dotsc$,
\begin{equation}\label{polygamma-sum}
\psi^{(n)}(z)=(-1)^{n+1}n!\sum_{k=0}^\infty\frac1{(z+k)^{n+1}}.
\end{equation}
\end{lem}

\section{Proofs of theorems and corollaries}

Now we are in a position to prove our theorems.

\begin{proof}[Proof of Theorem~\ref{batir-thm2.8-gen}]
Direct calculation gives
\begin{align}\label{phi'}
\phi'(x)&=\psi'(x)-\frac{e^{1/x}}{\bigl(e^{1/x}-1\bigr)x^2},\\
\phi''(x)&=\psi''(x)+\frac{2e^{1/x}} {\bigl(e^{1/x}-1\bigr)x^3} +\frac{e^{1/x}}{\bigl(e^{1/x}-1\bigr)x^4} -\frac{e^{2/x}}{\bigl(e^{1/x}-1\bigr)^2x^4} \label{phi''}
\end{align}
and
$$
\lim_{x\to\infty}\phi'(x)=\lim_{x\to\infty}\phi''(x)=0.
$$
\par
From the difference equation \eqref{psisymp4}, it is deduced that
\begin{equation}\label{phi-1-der}
\begin{split}
\phi'(x)-\phi'(x+1)&=\frac1{x^2}+ \frac1{\bigl[1-1/e^{1/(x+1)}\bigr](x+1)^2} -\frac1{\bigl(1-1/e^{1/x}\bigr)x^2}\\
&=\frac{e^{1/(x+1)}}{[e^{1/(x+1)}-1](x+1)^2}-\frac1{\bigl(e^{1/x}-1\bigr)x^2}.
\end{split}
\end{equation}
It is easy to see that $\phi'(x)-\phi'(x+1)>0$ for $x>0$ is equivalent to
\begin{equation}\label{equiv-ineq}
x^2\bigl(e^{1/x}-1\bigr)>(x+1)^2\bigl[1-e^{-1/(x+1)}\bigr].
\end{equation}
This can be expanded and simplified as
\begin{equation*}
\sum_{k=3}^\infty\frac1{k!}\biggl[\frac1{x^{k-2}} +\frac{(-1)^{k}}{(x+1)^{k-2}}\biggr]>0
\end{equation*}
which is valid clearly. By Lemma~\ref{key-lem}, it is concluded that the function $\phi'(x)$ is positive and $\phi(x)$ is strictly increasing on $(0,\infty)$.
\par
By utilization of \eqref{phi-1-der} and differentiation, it is acquired directly that
\begin{align*}
\phi''(x)-\phi''(x+1)&=[\phi'(x)-\phi'(x+1)]' \\* &=\frac{e^{1/(x+1)}\bigl[3+2x-2(x+1)e^{1/(x+1)}\bigr]} {\bigl[e^{1/(x+1)}-1\bigr]^2(x+1)^4} -\frac{(1-2x)e^{1/x}+2x}{\bigl(e^{1/x}-1\bigr)^2x^4}
\end{align*}
for $x>0$. It is obvious that the fact $\phi''(x)-\phi''(x+1)<0$ is equivalent to
\begin{equation*}
\biggl[\frac{e^{1/x}-1} {e^{1/(x+1)}-1}\biggr]^2 >\frac1{e^{1/(x+1)}}\biggl(\frac{x+1}x\biggr)^4 \frac{e^{1/x}-2x\bigl(e^{1/x}-1\bigr)} {1-2(1+x)\bigl[e^{1/(x+1)}-1\bigr]}.
\end{equation*}
Considering \eqref{equiv-ineq}, in order to prove above inequality, it is sufficient to show
\begin{equation*}
\frac1{e^{1/(x+1)}}>\frac{e^{1/x}-2x\bigl(e^{1/x}-1\bigr)} {1-2(1+x)\bigl[e^{1/(x+1)}-1\bigr]}
\end{equation*}
which is equivalent to
\begin{equation*}
3+2x-2e^{1/(x+1)}-(2x+1)e^{1/(x+1)+1/x}\triangleq h(x)<0.
\end{equation*}
Straightforward computation gives
\begin{gather*}
h'(x)=\frac{2x^2e^{1/(x+1)}+2(x+1)^2x^2 +\bigl(4x^2+4x+1-2x^4\bigr)e^{1/(x+1)+1/x}}{x^2(x+1)^2},\\
h''(x)=-\frac{2(2x+3)x^4+e^{1/x}\bigl(8x^5+26x^4+36x^3+24x^2+8x+1\bigr)} {x^4(x+1)^4e^{-1/(x+1)}}<0,
\end{gather*}
and $\lim_{x\to\infty}h'(x)=0$. Hence, the function $h'(x)$ for $x>0$ is decreasing and positive, and then the function $h(x)$ for $x>0$ is increasing. From $\lim_{x\to\infty}h(x)=-4$, it is deduced that $h(x)<-4<0$ for $x>0$. Consequently, utilizing Lemma~\ref{key-lem}, it is concluded that $\phi''(x)<0$ on $(0,\infty)$. The concavity of $\phi(x)$ is proved.
\par
From \eqref{binetcor}, it follows that
\begin{align*}
\phi(x)&=\ln x-\frac1x+\ln\bigl(e^{1/x}-1\bigr)+\int_{0}^{\infty} \left(\frac{1}{t}-\frac{1}{e^t-1}\right)e^{-xt}\td t\\
&=\ln\frac{\bigl(e^{1/x}-1\bigr)}{1/x}-\frac1x+\int_{0}^{\infty} \left(\frac{1}{t}-\frac{1}{e^t-1}\right)e^{-xt}\td t\\
&\to0
\end{align*}
as $x\to\infty$. Employing \eqref{psisymp4} for $i=1$ reveals
\begin{gather*}
\phi(x)=\psi(x)+\frac1x+\ln\bigl(e^{1/x}-1\bigr)-\frac1x =\psi(x+1)+\ln\bigl(e^{1/x}-1\bigr)-\ln e^{1/x}\\*
=\psi(x+1)+\ln\bigl(1-e^{-1/x}\bigr)\to\psi(1)=-\gamma
\end{gather*}
as $x\to0^+$. The proof of Theorem~\ref{batir-thm2.8-gen} is complete.
\end{proof}

\begin{proof}[Proof of Corollary~\ref{batir-thm2.8-gen-cor}]
These inequalities follow from the increasing monotonicity and concavity of $\phi(x)$ and formulas \eqref{phi'} and \eqref{phi''}.
\end{proof}

\begin{proof}[Proof of Theorem~\ref{batir-thm2.8-gen-2}]
Using \eqref{binetcor}, the function $\phi_\theta(x)$ becomes
\begin{align*}
\phi_\theta(x)&=\ln x-\frac1x+\ln\bigl(e^{\theta/x}-1\bigr)+\int_{0}^{\infty} \left(\frac{1}{t}-\frac{1}{e^t-1}\right)e^{-xt}\td t\\
&=\ln\frac{\bigl(e^{\theta/x}-1\bigr)}{1/x}-\frac1x+\int_{0}^{\infty} \left(\frac{1}{t}-\frac{1}{e^t-1}\right)e^{-xt}\td t\\*
&\to\ln\theta
\end{align*}
as $x\to\infty$. Employing \eqref{psisymp4} for $i=1$ reveals
\begin{align*}
\phi_\theta(x)&=\psi(x)+\frac1x+\ln\bigl(e^{\theta/x}-1\bigr)-\frac1x \\ &=\psi(x+1)+\ln\bigl(e^{\theta/x}-1\bigr)-\ln e^{1/x}\\
&=\psi(x+1)+\ln\bigl(e^{(\theta-1)/x}-e^{-1/x}\bigr)\\
&\to\begin{cases}\psi(1)
=-\gamma,&\theta=1\\
\infty,&\theta>1\\
-\infty,&0<\theta<1
\end{cases}
\end{align*}
as $x\to0^+$. The two limits in Theorem~\ref{batir-thm2.8-gen-2} is proved.
\par
Easy calculation yields
\begin{gather}
\label{phi''-c}
\phi'_\theta(x)=\psi'(x)-\frac{\theta e^{\theta/x}} {x^2\bigl(e^{\theta/x}-1\bigr)}
\triangleq\psi'(x)-\varphi(\theta,x),\\
\phi''_\theta(x)=\psi''(x)+\frac{\theta e^{\theta/x} \bigl[2x\bigl(e^{\theta/x}-1\bigr)-\theta\bigr]} {x^4\bigl(e^{\theta/x}-1\bigr)^2}
=\psi''(x) -\frac{\td\varphi(\theta,x)}{\td x}. \label{phi'''-c}
\end{gather}
\par
It is easy to verify that
\begin{equation*}
\frac{\td\varphi(\theta,x)}{\td \theta}=\frac{e^{\theta/x}\bigl(e^{\theta/x}-1 -\theta/x\bigr)}{x^2\bigl(e^{\theta/x}-1\bigr)^2}>0
\end{equation*}
and
\begin{align*}
\frac{\td{}^2\varphi(\theta,x)}{\td \theta\td x}
&=-\frac{e^{\theta/x} \bigl[(\theta/x)^2\bigl(e^{\theta/x}+1\bigr) -4(\theta/x)\bigl(e^{\theta/x}-1\bigr) +2\bigl(e^{\theta/x}-1\bigr)^2\bigr]} {x^3\bigl(e^{\theta/x}-1\bigr)^3}\\
&=-\frac{2\theta^2e^{\theta/x} \bigl\{\bigl(e^{\theta/x}-1\bigr)/2 +\bigl[\bigl(e^{\theta/x}-1\bigr)/(\theta/x)-1\bigr]^2\bigr\}} {x^5\bigl(e^{\theta/x}-1\bigr)^3}\\*
&<0
\end{align*}
for $\theta>0$ and $x>0$. These implies that the functions $\varphi(\theta,x)$ and $\frac{\td\varphi(\theta,x)}{\td x}$ are increasing and decreasing with $\theta>0$ respectively. Consequently, by using the two inequalities in Corollary~\ref{batir-thm2.8-gen-cor}, it is concluded for $0<\theta\le1$ and $x\in(0,\infty)$ that
\begin{equation*}
\varphi(\theta,x)\le\varphi(1,x) =\frac{e^{1/x}}{x^2\bigl(e^{1/x}-1\bigr)}<\psi'(x)
\end{equation*}
and
\begin{equation*}
\frac{\td\varphi(\theta,x)}{\td x}\ge\frac{\td\varphi(1,x)}{\td x} =\frac{e^{1/x}\bigl[1-2x\bigl(e^{1/x}-1\bigr)\bigr]} {\bigl(e^{1/x}-1\bigr)^2x^4}>\psi''(x).
\end{equation*}
Hence, the function $\phi'_\theta(x)$ is positive and $\phi''_\theta(x)<0$ on $(0,\infty)$ for $0<\theta\le1$, and the function $\phi_\theta(x)$ is increasing and concave on $(0,\infty)$ for $0<\theta\le1$.
\par
In order that $\phi'_\theta(x)<0$, by the right-hand side inequality in \eqref{psi'ineq}, it is sufficient to prove
\begin{equation*}
\frac1{x}+\frac1{x^2}-\frac{\theta e^{\theta/x}} {x^2\bigl(e^{\theta/x}-1\bigr)}\le0
\end{equation*}
which is equivalent to
$$
1+\frac{u}\theta-\frac{ue^{u}}{e^{u}-1}\le0
$$
for $u=\frac\theta{x}>0$. Therefore, it suffices to let
$$
\theta\ge\frac{u(e^u-1)}{1+(u-1)e^u}\triangleq\delta(u)
$$
for $u>0$. Since $\delta(u)$ is decreasing and $1<\delta(u)<2$ on $(0,\infty)$, when $\theta\ge2$, the function $\phi_\theta(x)$ is decreasing on $(0,\infty)$.
\par
In order that $\phi''_\theta(x)>0$, by the right-hand side inequality in \eqref{psi'ineq}, it is sufficient to show
\begin{equation*}
\frac{\theta e^{\theta/x} \bigl[2x\bigl(e^{\theta/x}-1\bigr)-\theta\bigr]} {x^4\bigl(e^{\theta/x}-1\bigr)^2}\ge\frac{1}{x^2}+\frac{2}{x^{3}}
\end{equation*}
which is equivalent to
\begin{equation*}
\frac{2ue^u(e^u-1-u/2)}{(e^u-1)^2}\ge1+\frac{2u}\theta
\end{equation*}
for $u=\frac\theta{x}>0$. Therefore, it suffices to let
\begin{equation*}
\frac2\theta\le\frac1u\biggl[\frac{2ue^u(e^u-1-u/2)}{(e^u-1)^2}-1\biggr] \triangleq\rho(u)
\end{equation*}
for $u>0$. Since $\rho(u)$ is increasing and $1<\rho(u)<2$ on $(0,\infty)$, it is sufficient to let $\theta\ge2$.
\par
If $\phi_\theta(x)$ is increasing on $(0,\infty)$, then $\phi'_\theta(x)>0$ means that
$$
[x^2\psi'(x)-\theta]e^{\theta/x}>x^2\psi'(x).
$$
It is well-known that $\psi'(x)>0$ on $(0,\infty)$, thus it is necessary that $\theta<x^2\psi'(x)$. Lemma~\ref{comp-thm-1} for $k=1$ gives
$$
\psi'(x)<\frac1{x}+\frac1{x^{2}}
$$
on $(0,\infty)$, and then $\theta<x+1$ on $(0,\infty)$. Hence, the required necessary condition $\theta\le1$ is proved.
\end{proof}

\begin{proof}[Proof of Corollary~\ref{batir-thm2.8-gen-cor-2}]
These inequalities follow directly from the monotonicity and convexity of $\phi_\theta(x)$ and formulas \eqref{phi''-c} and \eqref{phi'''-c}.
\end{proof}

\begin{proof}[Proof of Theorem~\ref{batir-thm-corr}]
By \eqref{psi-sum}, it follows that
\begin{equation}\label{batir-2.28}
\psi(x+1)=-\gamma+\sum_{k=1}^\infty\biggl(\frac1k-\frac1{k+x}\biggr)
\end{equation}
for $x>-1$. By the mean value theorem for differentiation, it is obvious that there exists a number $\mu=\mu(k)=\mu(k,x)$ for $k\in\mathbb{N}$ such that $-1<\mu(k)<x$ and
\begin{equation}\label{batir-2.29}
\frac1k-\frac1{k+x}=\frac{x}{[k+\mu(k)]^2}.
\end{equation}
Employing \eqref{batir-2.29} in \eqref{batir-2.28} leads to
\begin{equation}\label{batir-2.30}
\psi(x+1)=-\gamma+x\sum_{k=1}^\infty\frac1{[k+\mu(k)]^2}
\end{equation}
for $x>-1$. From \eqref{batir-2.29}, it is deduced that
$$
\mu(k)=\sqrt{k(k+x)}\,-k.
$$
It is not difficult to show that the mapping $k\to\mu(k)$ is strictly increasing on $[1,\infty)$ with
$$
\mu(1)=\sqrt{1+x}-1>-1\quad \text{and} \quad \lim_{k\to\infty}\mu(k)=\frac{x}2.
$$
Hence, from \eqref{batir-2.30} and \eqref{polygamma-sum}, it is concluded that inequality
\begin{multline*}
x\psi'\biggl(1+\frac{x}2\biggr) =x\psi'(1+\lim_{k\to\infty}\mu(k)) =x\sum_{k=1}^\infty\frac1{[k+\lim_{k\to\infty}\mu(k)]^2}\\ <\gamma+\psi(x+1) <x\sum_{k=1}^\infty\frac1{[k+\mu(1)]^2}=x\psi'(1+\mu(1)) =x\psi'\bigl(\sqrt{x+1}\,\bigr)
\end{multline*}
holds for $x>0$ and reverses for $-1<x<0$. The proof of Theorem~\ref{batir-thm-corr} is finished.
\end{proof}

\begin{proof}[Proof of Theorem~\ref{f-g-mon}]
Direct computation and utilization of the mean value theorem for differentiation gives
\begin{align*}
f'(x)&=\psi'(x+1)-\psi'\biggl(1+\frac{x}2\biggr) -\frac{x}2\psi''\biggl(1+\frac{x}2\biggr)\\
&=\frac{x}2\psi''(1+\xi(x))-\frac{x}2\psi''\biggl(1+\frac{x}2\biggr)\\
&=\frac{x}2\biggl[\psi''(1+\xi(x))-\psi''\biggl(1+\frac{x}2\biggr)\biggr],
\end{align*}
where $\xi(x)$ is between $\frac{x}2$ and $x$ for $x>-1$. Since $\psi''(x)$ is strictly increasing on $(0,\infty)$, it follows clearly that $f'(x)>0$ for $x\ne0$. Hence, the function $f(x)$ is strictly increasing on $(-1,\infty)$.
\par
Standard argument leads to
\begin{align*}
g'(x)&=\psi'\bigl(\sqrt{x+1}\,\bigr) -\psi'(x+1)+\frac{x}{2\sqrt{x+1}\,}\psi''\bigl(\sqrt{x+1}\,\bigr)\\
&=\frac{2u\bigl[\psi'(u)-\psi'\bigl(u^2\bigr)\bigr] +\bigl(u^2-1\bigr)\psi''(u)}{2u}\\
&=\frac{\bigl[\bigl(u^2-1\bigr)\psi'(u)-\psi\bigl(u^2\bigr)\bigl]'}{2u}
\end{align*}
for $x>-1$ and $u=\sqrt{x+1}\,>0$. Utilization of formulas \eqref{psi-sum} and \eqref{polygamma-sum} and direct differentiation gives
\begin{align}
\bigl(u^2-1\bigr)\psi'(u)-\psi\bigl(u^2\bigr)
&=\bigl(u^2-1\bigr) \sum_{i=0}^{\infty}\frac{1}{(u+i)^2}+\gamma -\sum_{i=0}^{\infty}\left(\frac{1}{1+i}-\frac{1}{u^2+i}\right)\notag\\
&=\gamma+\bigl(u^2-1\bigr)\sum_{i=0}^\infty\biggl[\frac{1}{(u+i)^2} -\frac1{(i+1)(u^2+i)}\biggr]\notag\\
&=\gamma+\sum_{i=1}^\infty\frac{i(u^2-1)(u-1)^2}{(i+1)(u+i)^2(u^2+i)} \label{u2-der}
\end{align}
and
\begin{equation}\label{u2-der-1}
\bigl[\bigl(u^2-1\bigr)\psi'(u)-\psi\bigl(u^2\bigr)\bigr]'
=\sum_{i=1}^\infty\frac{2i(u-1)^2 \bigl(u^3+2u^2+2iu+i\bigr)} {(i+u)^3\bigl(u^2+i\bigr)^2}>0
\end{equation}
for $u>0$. Hence, the function $g(x)$ is strictly increasing on $(-1,\infty)$.
\par
It is apparent that
$$
\lim_{x\to-1^+}f(x)=\lim_{x\to-1^+}\psi(x+1) +\psi'\biggl(\frac12\biggr)=-\infty.
$$
\par
By \eqref{psi'ineq} for $k=1$ and $\lim_{x\to\infty}\psi(x)=\infty$, it is clear that $\lim_{x\to\infty}f(x)=\infty$.
\par
By using \eqref{u2-der}, it is easy to see that
\begin{gather*}
\lim_{x\to-1^+}g(x)=\lim_{u\to0^+}\bigl[\bigl(u^2-1\bigr)\psi'(u) -\psi\bigl(u^2\bigr)\bigr]\\
=\gamma+\lim_{u\to0^+}\sum_{i=1}^\infty\frac{i(u^2-1)(u-1)^2} {(i+1)(u+i)^2(u^2+i)}
=\gamma-\sum_{i=1}^\infty\frac1{(i+1)i^2}
=1+\gamma-\frac{\pi^2}6
\end{gather*}
and
\begin{gather*}
\lim_{x\to\infty}g(x)=\lim_{u\to\infty}\bigl[\bigl(u^2-1\bigr)\psi'(u) -\psi\bigl(u^2\bigr)\bigr]\\
=\gamma+\lim_{u\to\infty}\sum_{i=1}^\infty\frac{i(u^2-1)(u-1)^2} {(i+1)(u+i)^2(u^2+i)}
=\gamma+\sum_{i=1}^\infty\frac{i}{i+1}
=\infty.
\end{gather*}
The proof of Theorem~\ref{f-g-mon} is complete.
\end{proof}

\begin{proof}[Proof of Theorem~\ref{correction-thm-3}]
It is not difficult to see that the factor $u^3+2u^2+2iu+i$ for $i\in\mathbb{N}$ in formula \eqref{u2-der-1} is positive if and only if $u>-1$. Therefore, the function $h(x)$ is strictly increasing on $(-1,\infty)$.
\par
The limits can be derived from \eqref{psi-sum} and \eqref{polygamma-sum}.
\end{proof}

\begin{proof}[Proof of Theorem~\ref{f-g-mon-poly}]
It is obvious that
\begin{align*}
[f_i(x)]'&=\psi^{(i+1)}(x+1)-\psi^{(i+1)}\biggl(1+\frac{x}2\biggr) -\frac{x}2\psi^{(i+2)}\biggl(1+\frac{x}2\biggr)\\
&=\frac{x}2\biggl[\psi^{(i+2)}(1+\eta(x)) -\psi^{(i+2)}\biggl(1+\frac{x}2\biggr)\biggr],
\end{align*}
where $\eta(x)$ is between $\frac{x}2$ and $x$. This means $(-1)^i[f_i(x)]'\ge0$ for $i\in\mathbb{N}$, and then the monotonicities of $f_i(x)$ on $(-1,\infty)$ for $i\in\mathbb{N}$ are proved.
\par
The two limits can be deduced easily from inequality \eqref{psi'ineq}.
\end{proof}


\begin{thebibliography}{99}

\bibitem{abram}
M. Abramowitz and I. A. Stegun (Eds), \textit{Handbook of Mathematical
Functions with Formulas, Graphs, and Mathematical Tables}, National Bureau of
Standards, Applied Mathematics Series \textbf{55}, 4th printing, with
corrections, Washington, 1965.

\bibitem{alzer-expo-math-2006}
H. Alzer, \textit{Sharp inequalities for the harmonic numbers}, Expo. Math. \textbf{24} (2006), no.~4, 385\nobreakdash--388.

\bibitem{Atanassov}
R. D. Atanassov and U. V. Tsoukrovski, \textit{Some properties of a class of
logarithmically completely monotonic functions}, C. R. Acad. Bulgare Sci.
\textbf{41} (1988), no.~2, 21\nobreakdash--23.

\bibitem{batir-jmaa-06-05-065}
N. Batir, \textit{On some properties of digamma and polygamma functions}, J. Math. Anal. Appl. \textbf{328} (2007), no.~1, 452\nobreakdash--465.

\bibitem{batir-new}
N. Batir, \textit{Some new inequalities for gamma and polygamma functions}, J.
Inequal. Pure Appl. Math. \textbf{6} (2005), no.~4, Art.~103; Available online
at \url{http://jipam.vu.edu.au/article.php?sid=577}.

\bibitem{CBerg}
C. Berg, \textit{Integral representation of some functions related to the
gamma function}, Mediterr. J. Math. \textbf{1} (2004), no.~4,
433\nobreakdash--439.

\bibitem{chen-qi-log-jmaa}
Ch.-P. Chen and F. Qi, \textit{Logarithmically completely monotonic functions
relating to the gamma function}, J. Math. Anal. Appl. \textbf{321} (2006), no.~1, 405\nobreakdash--411.

\bibitem{grin-ismail}
A. Z. Grinshpan and M. E. H. Ismail, \textit{Completely monotonic functions
involving the gamma and $q$\nobreakdash-gamma functions}, Proc. Amer. Math.
Soc. \textbf{134} (2006), 1153\nobreakdash--1160.

\bibitem{notes-best-simple-equiv.tex}
F. Qi, \textit{A completely monotonic function involving divided difference of
psi function and an equivalent inequality involving sums}, ANZIAM J. \textbf{48}
(2007), no.~4, 523\nobreakdash--532.

\bibitem{clark-ismail}
F. Qi, \textit{Certain logarithmically $N$-alternating monotonic functions involving gamma and $q$\nobreakdash-gamma functions}, Nonlinear Funct. Anal. Appl. \textbf{12} (2007), no.~4, 675\nobreakdash--685.

\bibitem{sandor-gamma-2}
F. Qi, \textit{Three classes of logarithmically completely monotonic functions involving gamma and psi functions}, Integral Transforms Spec. Funct. \textbf{18} (2007), no.~7, 503\nobreakdash--509.

\bibitem{compmon2}
F. Qi and Ch.-P. Chen, \textit{A complete monotonicity property of the gamma
function}, J. Math. Anal. Appl. \textbf{296} (2004), no.~2, 603\nobreakdash--607.

\bibitem{minus-one}
F. Qi and B.-N. Guo, \textit{Complete monotonicities of functions involving
the gamma and digamma functions}, RGMIA Res. Rep. Coll. \textbf{7} (2004),
no.~1, Art.~8, 63\nobreakdash--72; Available online at
\url{http://www.staff.vu.edu.au/rgmia/v7n1.asp}.

\bibitem{e-gam-rat-comp-mon}
F. Qi, B.-N. Guo and Ch.-P. Chen, \textit{Some completely monotonic functions
involving the gamma and polygamma functions}, J. Aust. Math. Soc.
\textbf{80} (2006), 81\nobreakdash--88.

\bibitem{auscmrgmia}
F. Qi, B.-N. Guo and Ch.-P. Chen, \textit{Some completely monotonic functions
involving the gamma and polygamma functions}, RGMIA Res. Rep. Coll. \textbf{7}
(2004), no.~1, Art.~5, 31\nobreakdash--36; Available online at
\url{http://www.staff.vu.edu.au/rgmia/v7n1.asp}.

\bibitem{widder}
D. V. Widder, \textit{The Laplace Transform}, Princeton University Press,
Princeton, 1941.

\end{thebibliography}
\end{document}